\documentclass[12pt]{article}
\usepackage{amsmath,amssymb}

\def\MC{\mathcal{M}}

\def\LC{\mathcal{L}}

\def\PC{\mathcal{P}}

\def\R{\mathbf{R}}

\def\E{\mathbf{E}}
\def\1{\mathbf{1}}

\newcommand{\Pf}{\par\noindent{\em Proof. }}

\def\pa{\partial}
\def\ep{\epsilon}
\def\tr{\rm{tr}}

\newtheorem{prop}{Proposition}[section]
\newtheorem{theorem}{Theorem}[section]

\begin{document}
\title{The L\'evy-Khintchine type operators with variable Lipschitz continuous coefficients
generate linear or nonlinear Markov processes and semigroups
 \thanks{submitted to PTRF September 2008}}
\author{Vassili N. Kolokoltsov\thanks{Department of Statistics, University of Warwick,
 Coventry CV4 7AL UK,
  Email: v.kolokoltsov@warwick.ac.uk}}
\maketitle

\begin{abstract} Ito's construction of Markovian solutions to stochastic equations driven by a L\'evy
noise is extended to nonlinear distribution dependent integrands
aiming at the effective construction of linear and nonlinear Markov
semigroups and the corresponding processes with a given
pseudo-differential generator. It is shown that a conditionally
positive integro-differential operator (of the L\'evy-Khintchine
type) with variable coefficients (diffusion, drift and L\'evy
measure) depending Lipschitz continuously on its parameters
(position and/or its distribution) generates a linear or nonlinear
Markov semigroup, where the measures are metricized by the
Wasserstein-Kantorovich metrics. This is a nontrivial but natural
extension to general Markov processes of a long known fact for
ordinary diffusions.
\end{abstract}

\paragraph{Key words.} Stochastic equations driven by L\'evy noise,
nonlinear integrators, Wasserstein-Kantorovich metric,
pseudo-differential generators, linear and nonlinear Markov
semigroups.

\paragraph{Running Head:} L\'evy-Khintchine operators with Lipschitz coefficients

\section{Introduction and formulation of main results}\label{sec1}

By $C(\R^n)$ (respectively $C_\infty(\R^n)$) we denote the Banach
space of continuous bounded functions on $\R^n$ (respectively its
subspace of functions vanishing at infinity) with the sup-norm
denoted by $\|\cdot\|$, and $C^k(\R^n)$ (resp. $C^k_c(\R^n)$)
denotes the Banach space of $k$ times continuously differentiable
functions with bounded derivatives on $\R^n$ (resp. its subspace of
functions with a compact support) with the norm being the sum of the
sup-norms of a function and all its partial derivative up to and
including the order $k$).

 For an $f\in C^1(\R^n)$ the gradient will be denoted by
 \[
\nabla f=(\nabla_1f,...,\nabla_nf)=(\frac{\pa f} {\pa
x_1},...,\frac{\pa f} {\pa x_1}).
\]
 For a measure $\nu $ and a mapping $F$ we denote by $\nu^F$ the
push forward of $\nu$ with respect to $F$ defined as $\nu^F(A)=\nu
(F^{-1}(A))$.

Further basic notations: ${\bf 1}_M$ is the indicator function of a
set $M$, $\MC (\R^d)$ is the set of finite positive Borel measures
on $\R^d$, $B_r$ is the ball of radius $r$ centered at the origin,
and the pairing $(f,\mu)$ for $f\in C(\R^d)$, $\mu \in \MC(\R^d)$
denotes the usual integration.
 The bold letters ${\bf E}$ and ${\bf P}$ will
denote expectation and probability. A positive number in the square
bracket, say $[x]$, will denote the integer part of it. By the small
letter $c$
 we shall denote various constants indicating in brackets
 (when appropriate) the parameters on which they depend.

 It is well known (the Courr\`ege theorem, see e.g. \cite{Ja}) that the generator $L$
 of a conservative (i.e. preserving constants) Feller
semigroup in $\R^d$ is conditionally positive ($f\ge 0, f(x)=0
\implies Lf(x)\ge 0$) and if its domain contains the space
$C_c^2(\R^d)$, then it has the following L\'evy-Khintchine form with
variable coefficients:
\begin{equation}
\label{fellergenerator}
 Lf(x)=\frac{1}{2}(G(x)\nabla,\nabla)f(x)+ (b(x),\nabla f(x))+\int
 (f(x+y)-f(x)-(\nabla f (x), y){\bf 1}_{B_1}(y))\nu (x,dy),
 \end{equation}
 where $G(x)$ is a symmetric non-negative matrix and $\nu(x,.)$ a Borel
 measure on $\R^d$ (called L\'evy measure) such that
\begin{equation}
\label{condlevy0}
 \int_{\R^n} \min (1,|y|^2)\nu (x;dy) <\infty, \quad \nu (\{0\})=0.
\end{equation}
The inverse question on whether a given operator of this form (or
better to say its closure) actually generates a Feller semigroup is
nontrivial and attracted lots of attention. One can distinguish
analytic and probabilistic approaches to this problem. The existence
results obtained by analytic techniques require certain
non-degeneracy condition on $\nu$, e.g. a lower bound for the symbol
of pseudo-differential operator $L$ (see e.g. \cite{Ba}, \cite{BC},
\cite{Ja}- \cite{Ko3} and references therein), and for the
construction of the processes via usual stochastic calculus one
needs to have a family of transformations $F_x$ of $\R^d$ preserving
the origin, regularly depending on $x$ and pushing a certain L\'evy
measure $\nu$ to the L\'evy measures $\nu(x,.)$, i.e.
$\nu(x,.)=\nu^{F_x}$ (see e.g. \cite{Ap}, \cite{BGJ}, \cite{Str}).
Of course yet more nontrivial is the problem of constructing the so
called nonlinear Markov semigroups solving the weak equations of the
form
\begin{equation}
\label{eqofnonlinearmarkovsemigroup}
 \frac{d}{dt}(f,\mu_t)=(L_{\mu_t}f,\mu_t), \quad \mu_t \in \PC
 (\R^d), \quad \mu_0=\mu,
\end{equation}
that should hold, say, for all $f\in C^2_c(\R^d)$, where $L_{\mu}$
has form \eqref{fellergenerator}, but with all coefficients
additionally depending on $\mu$, i.e.
\[
 L_{\mu}f(x)=\frac{1}{2}(G(x,\mu)\nabla,\nabla)f(x)+ (b(x,\mu),\nabla f(x))
 \]
\begin{equation}
\label{eqnonlinearMarkovgen}
 +\int
 (f(x+y)-f(x)-(\nabla f (x), y){\bf 1}_{B_1}(y))\nu (x,\mu,dy).
 \end{equation}
Equations of type \eqref{eqofnonlinearmarkovsemigroup} play
indispensable role in the theory of interacting particles (mean
field approximation) and exhaust all positivity preserving
evolutions on measures subject to certain mild regularity
assumptions (see e.g. \cite{Ko4}, \cite{Str}). A resolving semigroup
$U_t:\mu \mapsto \mu_t$ of the Cauchy problem for equation
\eqref{eqofnonlinearmarkovsemigroup} specified a so called
generalized or nonlinear Markov process $X(t)$, whose distribution
$\mu_t$ at time $t$ can be determined by the formula $U_{t-s}\mu_s$
from its distribution $\mu_s$ at any previous moment $s$.

 In the case of diffusions (when $\nu$ vanishes in \eqref{fellergenerator}
or \eqref{eqnonlinearMarkovgen}) the theory of the corresponding
semigroups is well developed, see \cite{McK} and more recent
achievements in \cite{GMN}. Also well developed is the case of pure
jump processes, see e.g. the treatment of the Boltzmann equation
(spatially trivial) in \cite{Szn}.

The goal of the present paper is to exploit the idea of nonlinear
integrators (see \cite{CN}, \cite{Ku}) combined with a certain
coupling of L\'evy processes in order to push forward the
probabilistic construction in a way that allows the natural
Lipschitz continuous dependence of the coefficients $G,b,\nu$ on
$x,\mu$ with measures equipped with their Wasserstein metric (see
the definition below). Thus obtained extension of the standard SDEs
with L\'evy noise represent a probabilistic counterpart of the
celebrated extension of the Monge mass transformation problem to the
generalized Kantorovich one. To streamline the exposition we shall
use Ito's approach (as exposed in detail in \cite{Str}) for
constructing the solutions of stochastic equations directly via
Euler approximation scheme bypassing the theory of stochastic
integration itself. Roughly speaking the idea is to approximate a
process with a given (formal) generator (or pre-generator) by
processes with piecewise L\'evy paths.

For a random variable $X$ we shall denote by $\LC (X)$ the
distribution (probability law) of $X$. Recall that the so called
Wasserstein-Kantorovich metrics $W_p$, $p\ge 1$, on the set of
probability measures $\PC (\R^d)$ on $\R^d$ are defined as
\begin{equation}
\label{eqdefmeasuredistance}
 W_p(\nu_1,\nu_2)=\left( \inf_{\nu} \int
|y_1-y_2|^p \nu (dy_1 dy_2)\right)^{1/p},
\end{equation}
where $\inf$ is taken over the class of probability measures $\nu$
on $\R^{2d}$ that couple $\nu_1$ and $\nu_2$, i.e. that satisfy
\begin{equation}
\label{eqcouplingcond}
 \int \int(\phi_1(y_1)+\phi_2(y_2))
\nu(dy_1dy_2)=(\phi_1,\nu_1)+(\phi_2,\nu_2)
\end{equation}
for all bounded measurable $\phi_1,\phi_2$.
 It follows directly from the definition that
\begin{equation}
\label{estimateonmeasuredistance}
 W_p^p (\mu,\mu') \le \E \|X-X'\|^p
\end{equation}
whenever  $\mu=\LC (X)$ and $\mu'=\LC (X')$.

 For random variable
$x,z$ we shall write sometimes shortly $W_p(x,z)$ for
$W_p(\LC(x),\LC(z))$ (with some obvious abuse of notation).

 It is well known (see e.g.
\cite{Vi}) that $(\PC (\R^d),W_p)$, $p\ge 1$ is a complete metric
space and that the convergence in this metric space is equivalent to
the convergence in the weak sense combined with the convergence of
the $p$th moments. In case $p=1$ the celebrated Monge-Kantorovich
theorem states that
\[
W_1(\mu_1,\mu_2)= \sup_{f\in Lip} |(f,\mu_1)-(f,\mu_2)|,
\]
where $Lip$ is the set of continuous functions $f$ such that
$|f(x)-f(y)|\le \|x-y\|$ for all $x,y$.

We shall need also the Wasserstein distances between the
distributions in the Skorohod space $D([0,T],\R^d)$ of cadlag paths
in $\R^d$ defined of course as
\begin{equation}
\label{eqdefmeasuredistanceonpaths}
 W_{p,T}(X_1,X_2)=\inf \left( \E \,
\sup_{t\le T}|X_1(t)-X_2(t)|^p \right)^{1/p},
\end{equation}
where $\inf$ is taken over all couplings of the distributions of the
random paths $X_1,X_2$. Notice that this distance is linked with the
uniform (and not Skorohod) topology on the path space.

To compare the L\'evy measures, we shall need an extension of these
distances to unbounded measures. Namely, let  $\MC_p (\R^d)$ denote
the class of Borel measures $\mu$ on $\R^d \setminus \{0\}$ (not
necessarily finite) with a finite $p$-th moment (i.e. such that
$\int |y|^p \mu (dy) <\infty$). For a pair of measures $\nu_1,\nu_2$
from  $\MC_p (\R^d)$ we define the distance $W_p(\nu_1,\nu_2)$ by
\eqref{eqdefmeasuredistance}, where $\inf$ is now taken over all
$\nu \in \MC_p (\R^{2d})$ such that \eqref{eqcouplingcond} holds for
all $\phi_1,\phi_2$ satisfying $\phi_i(.)/|.|^p \in C(\R^d)$. It is
easy to see that for finite measures this definition coincides with
the previous one and that if measures $\nu_1$ and $\nu_2$ are
infinite, the distance $W_p(\nu_1,\nu_2)$ is finite. \footnote{Let a
decreasing sequence of positive numbers $\ep_n^1$ be defined by the
condition that $\nu_1$ can be decomposed into the sum
$\nu_1=\sum_{n=1}^{\infty}\nu_1^n$ of the probability measures
$\nu_1^n$ having the support in the closed shells $\{x\in \R^d:
\ep_n^1 \le |x|\le \ep_{n-1}^1\}$ (where $\ep_0^1=\infty$).
Similarly $\ep_2^n$ and $\nu_2^n$ are defined. Then the sum
$\nu=\sum_{n=1}^{\infty}\nu_1^n\otimes \nu_2^n$ is a coupling of
$\nu_1$ and $\nu_2$ with a finite $\int |y_1-y_2|^p \nu(dy_1dy_2)$.}

 Moreover, by
the same argument as for finite measures (see \cite{RR} or
\cite{Vi}) one shows that whenever the distance $W_p(\nu_1,\nu_2)$
is finite, the infimum in \eqref{eqdefmeasuredistance} is achieved,
i.e. there exists a measure $\nu \in \MC_p (\R^{2d})$ such that

\begin{equation}
\label{eqmeasuredistanceattained}
 W_p(\mu_1,\mu_2)=\left( \int
|y_1-y_2|^p \nu (dy_1 dy_2)\right)^{1/p}.
\end{equation}

\begin{theorem}
\label{fellersemigroupslipschitzcoef}
 Let an operator $L$ have form \eqref{fellergenerator}, where
\begin{equation}
\label{eqlipschitzcoef}
  \|\sqrt {G(x_1)}-\sqrt {G(x_2)}\| + |b(x_1)-b(x_2)|+
  W_2(\1_{B_1}(.)\nu (x_1;.), \1_{B_1}(.)\nu(x_2;.))
  \le \kappa \|x_1-x_2\|
\end{equation}
with a certain constant $\kappa$, and
\begin{equation}
 \label{eqfamsymmsquareintegrableLevygenboundcoef0}
 \sup_x \left(\sqrt{G(x)}+ \vert b(x) \vert + \int_{B_1} |y|^2
  \nu(x, dy)\right) <\infty.
\end{equation}
Let the family of finite measures $\{\1_{\R^d\setminus B_1})(.)\nu
(x;.)\}$ be uniformly bounded, tight and depend weakly continuous on
$x$. Then $L$ extends to the generator of a conservative Feller
semigroup.
\end{theorem}

{\bf Remarks.} 1. The boundedness condition
\eqref{eqfamsymmsquareintegrableLevygenboundcoef0} is not essential
and can be dispensed with by the usual localization arguments, see
\cite{Ko6}. 2. Once the well posed-ness of the equations generated
by $L$ is obtained, it implies various extensions of the results on
the corresponding boundary value problems, problems with unbounded
coefficients, fractional dynamics or Malliavin calculus (see
\cite{Ta}, \cite{Ko3}, \cite{Ko5}, \cite{BGJ}) obtained earlier for
particular cases.

For example, assumption on $\nu$ is satisfied if one can decompose
the L\'evy measures $\nu(x;.)$ in the countable sums
$\nu(x;.)=\sum_{n=1}^\infty \nu_n(x;.)$ of probability measures so
that $W_2(\nu_i(x;.),\nu_i(z;.))\le a_i |x-z|$ and the series $\sum
a_i^2$ converges. It is well known that the optimal coupling of
probability measures (Kantorovich problem) can not always be
realized via a mass transportation (a solution to the Monge
problem), thus leading to the examples when the construction of the
process via standard stochastic calculus would not work. On the
other hand, no non-degeneracy is build in this example leading to
serious difficulties when trying to apply analytic techniques in
these circumstances.

Another important particular situation is that of a common star
shape of the measures $\nu(x;.)$, i.e. if they can be represented as
\begin{equation}
\label{eqstarshapes}
 \nu(x;dy)=\nu(x,s,dr)\, \omega (ds), \quad y\in \R^d, r=|y| \in \R_+,
 s=y/r \in S^{d-1},
\end{equation}
with a certain measure $\omega$ on $S^{d-1}$ and a family of
measures $\nu(x,s,dr)$ on $\R_+$. This allows to reduce the general
coupling problem to a much more easily handled one-dimensional one,
because evidently if $\nu_{x,y,s}(dr_1dr_2)$ is a coupling of
$\nu(x,s,dr)$ and $\nu(y,s,dr)$, then $\nu_{x,y,s}(dr_1dr_2)\omega
(ds)$ is a coupling of $\nu(x;.)$ and $\nu (y;.)$. If
one-dimensional measures have no atoms, their coupling can be
naturally organized via pushing along a certain mapping. Namely, the
measure $\nu^F$ is the pushing forward of a measure $\nu$ on $\R_+$
by a mapping $F:\R_+\mapsto \R_+$ whenever
\[
\int f(F(r))\nu (dr)=\int f(u) \nu^F(du)
\]
for a sufficiently rich class of test functions $f$, say for the
indicators of intervals. Suppose we are looking for a family of
monotone continuous bijections $F_{x,s}:\R_+\mapsto \R_+$ such that
$\nu^{F_{x,s}}=\nu(x,s,.)$. Choosing $f=\1_{[F(z),\infty)}$ as a
test function in the above definition of pushing yields
\begin{equation}
\label{eqdefpushinginonedim}
 G(x,s, F_{x,s}(z))=\nu ([z,\infty))
 \end{equation}
 for $G(x,s,z)= \nu (x,s,[z,\infty))=\int_z^{\infty}\nu (x,s,dy)$.
Clearly if all $\nu(x,s,.)$ and $\nu$ are unbounded, but bounded on
any interval separated from the origin, have no atoms and do not
vanish on any open interval, then this equation defines a unique
continuous monotone bijection $F_{x,s}:\R_+  \mapsto \R_+$ with also
continuous inverse. Hence we arrive to the following criterion.

 \begin{prop}
 \label{propcouplingforstarshapes}
Suppose the L\'evy measures $\nu(x;.)$ can be represented in the
form \eqref{eqstarshapes} and $\nu$ is a L\'evy measure on $\R_+$
such that  all $\nu(x,s,.)$ and $\nu$ are unbounded, have no atoms
 and do not vanish on any open interval. Then the family
 $\nu(x;.)$ depends Lipshitz continuous on $x$ in $W_2$ whenever the unique
 continuous solution $F_{x,s}(z)$ to \eqref{eqdefpushinginonedim} is
 Lipschitz continuous in $x$ with a constant $\kappa_F(z,s)$
 enjoying the condition
\begin{equation}
\label{eqlipforpushingF}
 \int_{\R_+} \int_{S^{d-1}} \kappa_F^2(r,s) \omega (ds) \nu (dr)
 <\infty.
 \end{equation}
\end{prop}

\Pf By the above discussion the solution $F$ specifies the coupling
$\nu_{x,y}(dr_1 dr_2 ds_1 ds_2)$ of $\nu(x;.)$ and $\nu(y;.)$ via
\[
\int f(r_1,r_2,s_1,s_2)\nu_{x,y}(dr_1 dr_2 ds_1 ds_2) =\int
f(F_{x,s}(r),F_{y,s}(r),s,s) \omega (ds) \nu (dr),
\]
so that for Lipschitz continuity of the family $\nu(x;.)$ it is
sufficient to have
\[
\int_{\R_+} \int_{S^{d-1}}
 (F_{x,s}-F_{y,s})^2 \omega (ds) \nu (dr)
\le c (x-y)^2,
\]
which is clearly satisfied whenever \eqref{eqlipforpushingF} holds.

The particular case of $\nu(x,s,.)$ above having densities with
respect to the Lebesgue measure on $\R_+$ is discussed in much
detail in \cite{Str}.

The point to make here is that a coupling for the sum of L\'evy
measures can be organized separately for each term allowing to use
the above statement for star shape components and, say, some
discrete methods for discrete parts.

Theorem \ref{fellersemigroupslipschitzcoef} is a straightforward
corollary of our main theorem that we shall formulate now. To make
our exposition more transparent we shall present the main arguments
in the case of $L_{\mu}$ having the form
\begin{equation}
\label{eqnonlinearMarkovgensecondmom}
  L_{\mu}f(x)=\frac{1}{2}(G(x,\mu)\nabla,\nabla ) f(x)+ (b(x,\mu),\nabla f(x))
 +\int
 (f(x+z)-f(x)-(\nabla f (x), z))\nu (x,\mu;dz)
 \end{equation}
with $\nu(x,\mu;.)\in
 \MC_2(\R^d)$. Let $Y_{\tau}(z,\mu)$ be a family of L\'evy processes
depending measurably on the points $z$ and probability measures
$\mu$ in $\R^d$ and specified by their generators
\[
L[z,\mu]f(x)=\frac{1}{2}(G(z,\mu)\nabla,\nabla)f(x)+
(b(z,\mu),\nabla f(x))
\]
\begin{equation}
\label{eqfamilyoflevygenerators}
 +\int (f(x+y)-f(x)-(\nabla f (x), y))\nu (z,\mu;dy)
 \end{equation}
 where $\nu (z,\mu) \in \MC_2(\R^d)$. Under the conditions of Theorem
 \ref{thnonlinearlevybasedstocheq} given below, the existence of such a
 family follows from the well known randomization
lemma \footnote{It states that if $\mu (x, dz)$ is a probability
kernel from a measurable space $X$ to a Borel space $Z$, then there
exists a measurable function $f:X\times [0,1] \to Z$ such that if
$\theta$ is uniformly distributed on $[0,1]$, then $f(X,\theta)$ has
distribution $\mu(x,.)$ for every $x\in X$.} (see e.g. \cite{Kal},
Lemma 3.22), because  by Proposition \ref{proplevycoupling} the
mapping from $z$, $\mu$ to the law of the L\'evy process
$Y_{\tau}(z,\mu)$ is continuous, hence measurable, and consequently,
by this Lemma (with $Z$ being the complete metric space
$D(\R_+,\R^d)$, and hence a Borel space) one can define all
$Y_{\tau}(z,\mu)$ on the single standard probability space $[0,1]$.
Let us stress for clarity that the processes $Y_{\tau}(x,\mu)$
depend on $x,\mu$ only via the parameters of the generator, i.e.,
say, the random variable $\xi=x+Y_{\tau}(x,\LC(x))$ has the
characteristic function
\[
\E e^{ip\xi}= \int \E e^{ip(x+Y_{\tau}(x,\LC(x))} \mu (dx).
\]

Our approach to solving \eqref{eqofnonlinearmarkovsemigroup} is via
the solution to the following nonlinear distribution dependent
stochastic equation with nonlinear L\'evy type integrators:
\begin{equation}
\label{stocheqfulnonlin}
 X(t)=X+\int_0^t dY_s (X(s),\LC(X(s))),\quad
\LC(X)=\mu,
\end{equation}
with a given initial distribution $\mu$ and a random variable $X$
independent of $Y_{\tau}(z,\mu)$.

We shall define the solution through the Euler type approximation
scheme, i.e. by means of the approximations $X_{\mu}^{\tau}$:
\begin{equation}
\label{eqdefapproximproc}
 X_{\mu}^{\tau}(t)
 = X_{\mu}^{\tau}(l\tau)+ Y_{t-l\tau}^l(
 X_{\mu}^{\tau}(l\tau), \LC (X_{\mu}^{\tau}(l\tau)) ), \quad \LC
 (X_{\mu}^{\tau}(0))=\mu,
\end{equation}
where $l\tau <t\le (l+1)\tau$, $l=0,1,2,...$, and 
$Y_{\tau}^l(x,\mu)$ is a collection (depending on $l$) of independent
families of the L\'evy processes $Y_{\tau}(x,\mu)$ introduced above.
Clearly these approximation processes are cadlag.

For $x\in \R^d$ we shall write shortly $X_x^{\tau}(k\tau)$ for
$X_{\delta_x}^{\tau}(k\tau)$.

By the weak solution to \eqref{stocheqfulnonlin} we shall mean the
weak limit of $X_{\mu}^{\tau_k}$, $\tau_k=2^{-k}$, $k\to \infty$, in
the sense of the distributions on the Skorohod space of cadlag paths
(which is of course implied by the convergence of the distributions
in the sense of the distance \eqref{eqdefmeasuredistanceonpaths}).
Alternatively one could define it as a solution to the corresponding
nonlinear martingale problem (see below the proof of the main
theorem) or directly via the construction of the corresponding
stochastic integral. This issue is addressed in detail in
\cite{Ko6}, our purpose here being the construction of a Markov
process with a given generator.

The following is our main result.

\begin{theorem}
\label{thnonlinearlevybasedstocheq}
 Let an operator $L_{\mu}$ have form
 \eqref{eqnonlinearMarkovgensecondmom}. Moreover
\begin{equation}
\label{eqlipschitzcoefnonlinearinW2}
  \|\sqrt {G(x,\mu)}-\sqrt {G(z,\eta)}\|+|b(x,\mu)-b(z,\eta)|+
  W_2(\nu (x,\mu;.), \nu(z,\eta;.))
  \le \kappa (|x-z|+W_2(\mu,\eta)),
\end{equation}
holds true with a constant $\kappa$ and
\begin{equation}
 \label{eqfamsymmsquareintegrableLevygenboundcoef}
 \sup_{x,\mu} \left(\sqrt{G(x,\mu)}+ \vert b(x,\mu) \vert + \int |y|^2
  \nu(x,\mu, dy)\right) <\infty.
\end{equation}
Then

(i) for any $\mu \in \PC(\R^d)\cap \MC_2(\R^d)$ the approximations
$X_{\mu}^{\tau_k}$ converge to a process $X_{\mu}(t)$ in the sense
that
\begin{equation}
\label{eqconvergeofapprinW2}
  \sup_{\mu} \sup_{t\in [0,t_0]} W^2_2\left( X_{\mu}^{\tau_k}([t/\tau_k]\tau_k, X_{\mu}(t)\right)
 \le c(t_0) \tau_k
\end{equation}
for any $t_0$, and even stronger
\begin{equation}
\label{eqconvergeofapprinW2onpaths}
  \sup_{\mu} W^2_{2,t_0}\left( X_{\mu}^{\tau_k}, X_{\mu}\right)
 \le c(t_0) \tau_k;
\end{equation}

(ii) the distributions $\mu_t=\LC (X_{\mu}(t))$ depend
$1/2$-H\"older continuous on $t$ in the metric $W_2$ and
$X_{\mu}(t)$ depend Lipschitz continuously on the initial condition
in the following sense:
 \begin{equation}
\label{eqcontinuityofsolutionprocinW2}
 \sup_{t\in [0,t_0]}
W_2^2(X_{\mu}(t),X_{\eta}(t)) \le c(t_0)W^2_2(\mu,\eta);
 \end{equation}

(iii) the processes
 \begin{equation}
 \label{eqnonlinearmartingales}
 M(t)=f(X_{\mu}(t))-f(X_{\mu}^0)
  -\int_0^{t}(L_{\LC (X_{\mu}(s))}f(X_{\mu}(s))
 \, ds
 \end{equation}
 are martingales for any $f\in C^2(\R^d)$; in other words, the
 process $X_{\mu}(t)$ solves the corresponding (nonlinear)
 martingale problem;

(iv) the distributions $\mu_t=\LC (X_{\mu}(t))$ satisfy the weak
nonlinear equation \eqref{eqofnonlinearmarkovsemigroup} (that holds
for all $f\in C^2(\R^d)$);

(v) the resolving operators $U_t: \mu \mapsto \mu_t$ of the Cauchy
problem \eqref{eqofnonlinearmarkovsemigroup} form a nonlinear Markov
semigroup, i.e. they are continuous mappings from $\PC (\R^d)\cap
\MC_2(\R^d)$ (equipped with the metric $W_2$) to itself such that
$U_0$ is the identity mapping and $U_{t+s}=U_tU_s$ for all $s,t\ge
0$. If $L[z,\mu]$ do not depend explicitly on $\mu$ the operators
$T_tf(x)=\E f(X_x(t))$ form a conservative Feller semigroup
preserving the space of Lipschitz continuous functions.

\end{theorem}

This theorem is proved in the next section. In Sections 3 we obtain
some regularity criteria for the Markov semigroups constructed.

A simple meaningful example is given by the nonlinear kinetic
equations
\begin{equation}
\label{meanfieldforpotentialinteration}
 \frac{d}{dt}(f,\mu_t)=(Lf,\mu_t)+\int (K(x,y),\nabla f (x))\mu_t
 (dx)\mu_t(dy),
 \end{equation}
with $L$ being of form \eqref{fellergenerator} with Lipschitz
continuous coefficients and $K$ being a bounded Lipschitz continuous
mapping $R^{2d} \mapsto R^d$, which arise as the mean-field limit
for potentially interacting Feller processes.

Theorem \ref{fellersemigroupslipschitzcoef} follows now from Theorem
\ref{thnonlinearlevybasedstocheq} by the standard perturbation
theory, since dividing the generator into two parts, where the first
part is the integral term with the L\'evy measure reduced to
$\R^d\setminus B_1$, one gets a sum of two generators, one of which
is bounded in $C_{\infty}(\R^d)$ (as follows from the assumed
tightness) and the other satisfies Theorem
\ref{thnonlinearlevybasedstocheq}.

It is worth noting that in a simpler case of generators of up to the
first order the continuity of L\'evy measures with respect to a more
easy handled metric $W_1$ is sufficient, as shows the following
result, whose proof is omitted (as being a simplified version of the
proof of Theorem \ref{thnonlinearlevybasedstocheq}).

\begin{theorem}
\label{thnonlinearlevybasedstocheq1}
 Let an operator $L_{\mu}$ have the form
\begin{equation}
\label{eqnonlinearMarkovgenfirstmom}
 L_{\mu}f(x)= (b(x,\mu),\nabla f(x))
 +\int
 (f(x+z)-f(x))\nu (x,\mu;dz), \quad \nu(x,\mu;.)\in
 \MC_1(\R^d).
 \end{equation}
and
\begin{equation}
\label{eqlipschitzcoefnonlinearinW1} \|b(x,\mu)-b(z,\eta)\|+
  W_1(\nu (x,\mu;.), \nu(z,\eta;.))
  \le \kappa (\|x-z\|+W_1(\mu,\eta))
\end{equation}
holds true with a constant $\kappa$. Then for any $\mu \in
\PC(\R^d)\cap \MC_1(\R^d)$ there exists a process $X_{\mu}(t)$
solving \eqref{stocheqfulnonlin} (with analogously defined
$Y_{\tau}(z,\mu)$) such that
\begin{equation}
 \label{eqconvergeofapprinW1onpaths}
  \sup_{\mu} W_{1,t_0}\left( X_{\mu}^{\tau_k}, X_{\mu}\right)
 \le c(t_0) \tau_k,
\end{equation}
 the distributions $\mu_t=\LC (X(t))$ depend $1/2$-H\"older
continuous on $t$ in the metric $W_1$ and $X_{\mu}(t)$ depend
Lipschitz continuously on the initial condition in the following
sense:
 \begin{equation}
 \label{eqcontinuityofsolutionprocinW1}
W_1(X_{\mu}(t),X_{\eta}(t)) \le c(t_0)W_1(\mu,\eta).
 \end{equation}
Moreover, the processes \eqref{eqnonlinearmartingales} are
martingales for any $f\in C^1(\R^d)$ and the distributions
$\mu_t=\LC (X_{\mu}(t))$ satisfy the weak nonlinear equation
\eqref{eqofnonlinearmarkovsemigroup} (that holds for all $f\in
C^1(\R^d)$). If $L[z,\mu]$ do not depend explicitly on $\mu$ the
operators $T_tf(x)=\E f(X_x(t))$ form a conservative Feller
semigroup.
\end{theorem}

In Appendix we describe a coupling of L\'evy processes that is
crucial for our purposes.

\section{Proof of Theorem \ref{thnonlinearlevybasedstocheq}}\label{sec2}

  {\bf Step 1} ({\it
uniform continuity of the approximations with respect to initial
data}).

One has
\[
W^2_2(x_1+Y_s(x_1,\LC(x_1)),x_2+Y_s(x_2,\LC(x_2)))\le
\E(\xi_1-\xi_2)^2
\]
for any random variable $(\xi_1,\xi_2)$ with the projections
$\xi_i=x_i+Y_s(x_i,\mu_i)$, $\mu_i=\LC(x_i)$, $i=1,2$. Let us choose the coupling
described by the characteristic function
\[
\E e^{i(p_1\xi_1+p_2\xi_2)} = \int_{\R^{4d}}
e^{ip_1(x_1+y_1)+ip_2(x_2+y_2)}\mu(dx_1dx_2)
P^s_{x_1,x_2,\mu_1,\mu_2}(dy_1dy_2),
\]
where $\mu$ is an arbitrary coupling of the
random variables $x_1,x_2$ and $P^s$ is the coupling of the
L\'evy processes $Y_s(x_i,\mu_i)$ given by Proposition
\ref{proplevycoupling}. Consequently,
\[
\E(\xi_1-\xi_2)^2=-\frac{d^2}{dp^2}\mid_{p=0} \E e^{ip(\xi_1-\xi_2)}
\]
\[
=\int_{\R^{4d}} [(x_1+y_1)-(x_2+y_2)]^2\mu(dx_1dx_2)
P^s_{x_1,x_2,\mu_1,\mu_2}(dy_1dy_2),
\]
which by \eqref{levycouplingestimate2} does not exceed
\[
\int_{\R^{2d}}
\left((x_1-x_2)^2+cs[(x_1-x_2)^2+W^2_2(\LC(x_1),\LC(x_2))]\right)\mu(dx_1dx_2).
\]
Consequently, by \eqref{estimateonmeasuredistance},
\begin{equation}
 \label{eqcontofonestep0}
 \E(\xi_1-\xi_2)^2 \le
\int_{\R^{2d}} (1+2cs)(x_1-x_2)^2\mu(dx_1dx_2).
\end{equation}
 Hence, taking infimum over all couplings, yields
\begin{equation}
 \label{eqcontofonestep}
W^2_2(x_1+Y_s(x_1,\LC(x_1)),x_2+Y_s(x_2,\LC(x_2)))\le
(1+2cs)W_2^2(\LC(x_1),\LC(x_2)).
\end{equation}
Applying this inequality inductively, yields
 \begin{equation}
 \label{equniformcontinuityofappr}
W_2^2(X^{\tau}_{\mu}(s),X^{\tau}_{\eta}(s)) \le e^{1+2c
s} W^2_2 (\mu,\eta)
 \end{equation}
 with a constant $c$ uniformly for all $\tau \le 1, s>0$, $\mu,\eta \in \PC(\R^d)\cap \MC_2(\R^d)$.  

 {\bf Step 2} ({\it subdivision and the existence of the limit}).

We want to estimate the $W_2$ distance between the random variables
\[
\xi_1=x+Y_{\tau}(x,\mu)=x'+ Y'_{\tau/2}(x,\mu), \quad
\xi_2=z'+Y'_{\tau/2}(z',\eta'),
\]
where the families $Y_s$ and $Y'_s$ are independent,
\[
x'=x+ Y_{\tau/2}(x,\mu), \quad z'=z+ Y_{\tau/2}(z,\eta),
\]
and $\mu=\LC (x)$,  $\eta=\LC (z)$, $\eta'=\LC (z')$. We shall
couple $\xi_1$ and $\xi_2$ using sequentially Proposition
\ref{proplevycoupling}. Namely, we shall define it by the equation
\[
\E f(\xi_1,\xi_2) =\int_{\R^{6d}}f(x+v_1+y_1, z+v_2+y_2) \mu (dx dz)
P^{\tau/2}_{x,z,\mu,\eta}(dv_1dv_2)P^{\tau/2}_{x,z',\mu,\eta'}
(dy_1dy_2)
\]
for $f\in C(\R^{2d})$, where, say, $P^{\tau/2}_{x,z',\mu,\eta'}$ is
the coupling of the L\'evy processes $Y'_{\tau/2}(x,\mu)$ and
$Y'_{\tau/2}(z',\eta')$ with $z'=z+v_2$ given by Proposition \ref{proplevycoupling}
(note that the probability law $\eta'$ is the function of $z,\eta$).

Now by \eqref{levycouplingestimate2}
\[
W_2^2(\xi_1,\xi_2) \le \E (\xi_1-\xi_2)^2
\]
\[
\le \E (x'-z')^2
 + c\tau [\E(x'-z')^2 +\E (x-z')^2+W_2^2(\mu,\eta')].
\]
Hence, by \eqref{eqcontofonestep0} and
\eqref{estimateonmeasuredistance} $W_2^2(\xi_1,\xi_2)$ does not
exceed
\[
W_2^2(x,z)(1+2c\tau)(1+c\tau)+2c\tau \E (x-z')^2
\]
and consequently also
 \[
 W_2^2(x,z)(1+2c\tau)+4c\tau \E(Y_{\tau/2}(z,\eta))^2
\]
(with another constant $c$) so that
 \[
  W_2^2(\xi_1,\xi_2)\le W_2^2(x,z)(1+c\tau)+c\tau^2
\]
(with yet another $c$), because the second moments of our processes
$Y_{\tau}$ are bounded due to assumption
\eqref{eqfamsymmsquareintegrableLevygenboundcoef}. Consequently
\begin{equation}
\label{eqestimatesubdivision}
 W_2^2(X^{\tau}_{\mu}(k\tau),X^{\tau
/2}_{\mu}(k\tau))
  \le c\tau^2 +(1+c\tau)W_2^2(X^{\tau}_{\mu}((k-1)\tau),X^{\tau
  /2}_{\mu}((k-1)\tau)).
\end{equation}
 By induction one estimates the l.h.s. of this inequality by
\[
\tau^2[1+(1+c\tau)+(1+c\tau)^2+...+(1+c\tau)^{(k-1)}] \le c^{-1}\tau
(1+c\tau)^k \le c(t_0) \tau.
\]
 Repeating this subdivision and using the triangle inequality for distances yields
\[
W_2^2(X^{\tau}_{\mu}(k\tau),X^{\tau /2^m}_{\mu}(k\tau))
  \le c(t_0)\tau.
\]
This implies the existence of the limit
$X_x^{\tau_k}([t/\tau_k]\tau_k)$, as $k\to \infty$, in the sense of
\eqref{eqconvergeofapprinW2}.

Observe now that \eqref{equniformcontinuityofappr} implies
\eqref{eqcontinuityofsolutionprocinW2}. Moreover, the mapping
$T_tf(x)=\E f(X_x(t))$ preserves the set of Lipschitz continuous
functions. In fact, if $f$ is Lipschitz with the constant $h$, then
\[
 |\E f(X_x^{\tau}([t/\tau]\tau))-\E f(X_z^{\tau}([t/\tau]\tau))|
 \le h\E\|X_x^{\tau}([t/\tau]\tau)-(X_z^{\tau}([t/\tau]\tau)\|
 \]
 \[
 \le h \left(\E\|X_x^{\tau}([t/\tau]\tau)-(X_z^{\tau}([t/\tau]\tau)\|^2
 \right)^{1/2}.
\]
for any coupling of the processes $X_x^{\tau}$ and $X_z^{\tau}$.
Hence by \eqref{equniformcontinuityofappr}
\[
 |\E f(X_x^{\tau}([t/\tau]\tau))-\E f(X_z^{\tau}([t/\tau]\tau))|
  \le hc(t_0) W_2(x,z).
 \]
In particular, $T_t$ preserves constant functions.
 Similarly one shows (first for Lipschitz continuous $f$ and then
 for all $f\in C_{\infty}(\R^d)$ via standard approximation) that
\begin{equation}
\label{th1.3}
 \sup_{t\in [0,t_0]}\sup_x |\E
f(X_x^{\tau_k}([t/\tau_k]\tau_k))-\E f(X_x(t))| \to 0, \quad k \to
\infty,
\end{equation}
for all $f\in C_{\infty}(\R^d)$. Moreover, as the dynamics of
averages of the approximation processes clearly preserve the space
$C_{\infty}(\R^d)$, the same holds for the limiting mappings $T_t$..
Consequently $T_tf=\E f(X_x(t))$ is a positivity preserving family
of contractions in $C(\R^d)$ that preserve constants and the space
$C_{\infty}(\R^d)$. Hence the mappings $U_t: \mu \mapsto \mu_t$ form
a (nonlinear) Markov semigroup, and if $L[z,\mu]$ do not depend
explicitly on $\mu$, the operators $T_tf(x)=\E f(X_x(t))$ form a
conservative Feller semigroup. The Markov (or semigroup) property of
the solutions follows from the construction (a detailed discussion
of this fact in a similar situation is given in \cite{Str}).

From the inequality
 \[
 W_2^2(\LC (X_{\mu}^{\tau}(l\tau)),\LC (X_{\mu}^{\tau}((l-1)\tau)))
 \le
 \E \left[ Y_{\tau}^{l-1}(X_{\mu}^{\tau}((l-1)\tau),
 \LC (X_{\mu}^{\tau}((l-1)\tau)))\right]^2\le c \tau
 \]
 it follows that the curve $\mu_t$ depends $1/2$-H\"older
 continuously on $t$ in $W_2$.

{\bf Step 3} ({\it improving convergence and solving the martingale
problem})

 The processes
\begin{equation}
 \label{eqapprmart}
 M_{\tau}(t)=f(X_{\mu}^{\tau}(t))-f(X)
  -\int_0^tL[X^{\tau}_{\mu}([s/\tau]\tau),\mu^{\tau}_{[s/\tau]}]
  f(X^{\tau}_{\mu}(s))
 \, ds, \quad \mu= \LC (X),
 \end{equation}
 where  $\mu_l^{\tau}=\LC (X^{\tau}_{\mu} (l\tau))$,
 are martingales by Dynkin's formula, applied to L\'evy processes
 $Y_{\tau}(z,\mu)$.
 Our aim is to pass to the limit
  $\tau_k \to 0$ to obtain the martingale characterization of the
  limiting process. But let us first strengthen our convergence
  result.

Observe that the step by step inductive coupling of the trajectories
$X^{\tau}_{\mu}$ and $X^{\tau}_{\eta}$ used above to prove
\eqref{equniformcontinuityofappr}  actually defines the coupling
between the distributions of these random trajectories in the
Skorohod space $D([0,t_0], \R^d)$ for any $t_0$, i.e. a random
trajectory $(X^{\tau}_{\mu},X^{\tau}_{\eta})$ in $D([0,t_0],
\R^{2d})$. One can construct the Dynkin martingales for this coupled
process in the same way as above for $X^{\tau}_\mu$. Namely, for a
function $f$ of two variables with bounded second derivatives the
process
\[
 M_{\tau}(t)=f(X_{\mu}^{\tau}(t),X_{\eta}^{\tau}(t))
  -\int_0^t\tilde L_sf(X^{\tau}_{\mu}(s), X^{\tau}_{\eta}(s))
 \, ds, \quad \mu= \LC (x_{\mu}), \, \eta=\LC(x_{\eta}),
 \]
is a martingale, where $\tilde L_t$ is the coupling operator
\eqref{levycouplinggenerator2} constructed from the L\'evy processes
$Y$ with parameters
$X^{\tau}_{\mu}([t/\tau]\tau),\mu^{\tau}_{[t/\tau]}$ and
$X^{\tau}_{\eta}([t/\tau]\tau),\eta^{\tau}_{[t/\tau]}$.

Choosing
$f(x,y)=(x-y)^2$ leads to the martingale of the form
\[
(X_{\mu}^{\tau}(t)-X_{\nu}^{\tau}(t))^2
 +\int_0^t O(1)(X_{\mu}^{\tau}(s)-X_{\nu}^{\tau}(s))^2 \, ds
  \]
  (the estimate for the integrand follows from \eqref{levycouplinggenerator2b} and the assumed Lipschitz continuity of the coefficients of $L$). 
Applying the martingale property in conjunction with Gronwall's
lemma yields
\begin{equation}
 \label{equniformcontinuityofappr1}
\sup_{s\le t} \E (X_{\mu}^{\tau}(s)-X_{\nu}^{\tau}(s))^2 \le c(t) \E
(X_{\mu}(0)-X_{\nu}(0))^2,
 \end{equation}
giving another proof of \eqref{equniformcontinuityofappr}. Moreover,
applying Doob's maximal inequality (with $p=2$) to the vector-valued martingale
of the form
\[
\tilde M_{\tau}(t)=X_{\mu}^{\tau}(t)-X_{\eta}^{\tau}(t)
 +\int_0^t O(1)\vert X_{\mu}^{\tau}(s)-X_{\eta}^{\tau}(s)\vert \, ds
  \]
 constructed from $f(x,y)=x-y$ and using \eqref{equniformcontinuityofappr1} yields 
\[
\E \sup_{s\le t} |\tilde M_{\tau}(s)|^2 \le c(t) \E
(X_{\mu}(0)-X_{\nu}(0))^2,
\]
which in turn implies
\[
\E \sup_{s\le t} (X^{\tau}_{\mu}(s)-X^{\tau}_{\eta}(s))^2 \le c(t) \E
(X_{\mu}(0)-X_{\nu}(0))^2.
 \]
This allows to improve \eqref{equniformcontinuityofappr} to the
estimate of the distance on paths:
\begin{equation}
 \label{equniformcontinuityofappronpaths}
W_{2,T}^2(X^{\tau}_{\mu},X^{\tau}_{\eta})^2 \le c(T) W^2_2
(\mu,\eta).
 \end{equation}
 Similarly one can strengthen the estimates for subdivisions leading
 to the convergence of the distributions on paths
 \eqref{eqconvergeofapprinW2onpaths}.

Using the Skorohod theorem for the weak converging sequence of
random trajectories $X_{\mu}^{\tau_k}$ (let us stress again that the
convergence with respect to the distance
\eqref{eqdefmeasuredistanceonpaths} implies the weak convergence of
the distributions in the sense of the Skorohod topology), one can
put them all on a single probability space forcing the processes
$X^{\tau_k}_{\mu}$ to converge to $X_{\mu}$ almost surely in the
sense of the Skorohod topology.

Passing to the limit $\tau=\tau_k\to 0$ in \eqref{eqapprmart}, using
the continuity and boundedness of $f$ and $Lf$ and the dominated
convergence theorem  allows to conclude that the martingales
$M_{\tau}(t)$ converge almost surely and in $L^1$ to the martingale
 \[
 M(t)=f(X_{\mu}(t))-f(X)
  -\int_0^{t}(L_{\LC (X_{\mu}(s))}f)(X_{\mu}(s))
 \, ds,
 \]
 in other words that the process $X_{\mu}(t)$ solves the corresponding
 (nonlinear) martingale problem.

 {\bf Step 4} ({\it completion})

 To prove \eqref{eqofnonlinearmarkovsemigroup} one writes using the martingale
 properties of $M(t)$:
 \[
 \frac{d}{dt} (f,\mu_t)=\lim_{s\to 0}
 \frac{1}{s}\E (f,X_{\mu}(t+s)-X_{\mu}(t))=
\lim_{s\to 0}
 \frac{1}{s} \E \int_t^{t+s}(L_{\LC (X_{\mu}(s))}f)(X_{\mu}(s))
 \, ds
 \]
 \[
 =(L_{\mu_t}f, \mu_t))+\lim_{s\to 0}
 \frac{1}{s} \E \int_t^{t+s}[(L_{\LC (X_{\mu}(s))}f)(X_{\mu}(s))
 -L_{\mu_t}f(X_{\mu}(t)]
 \, ds,
 \]
 implying \eqref{eqofnonlinearmarkovsemigroup} by the continuity of $\mu_t$.

\section{Regularity}

Discussing regularity we reduce our attention for simplicity to
Feller processes. It is known (see e.g. \cite{Ko4}) that from the
sufficient regularity of nonhomogeneous versions of these Feller processes one
can naturally deduce the uniqueness and regularity for the
corresponding nonlinear problems.

 By $C^k_{Lip}$ (respectively
$C^k_{\infty}$) we shall denote the subspace of functions from
$C^k(\R^d)$ with a Lipschitz continuous derivative of order $k$
(respectively with all derivatives up to order $k$ vanishing at
infinity).

We shall discuss in detail only the first derivative.

\begin{theorem}
\label{thstocheqreg1stder} Assume the conditions of Propositions
\ref{proplevyderivcont} and \ref{proplevyderivcont1} hold.  Then the
spaces $C^1_{Lip}$ and $C^1_{Lip}\cap C^1_{\infty}$ are invariant
under the semigroup $T_t$ constructed above from the generator
\begin{equation}
\label{eqMarkovgenwithsecondmom}
Lf(x)=\frac{1}{2}(G(x)\nabla,\nabla)f(x)+ (b(x),\nabla f(x))
 +\int (f(x+y)-f(x)-(\nabla f(x), y))\nu (x,dy),
 \end{equation}
  and for any $f\in C^1_{Lip}$,
$\phi \in L_1\cap C_{\infty}(\R^d)$
\begin{equation}
\label{eqevolutioneqforsemigroupinLlip}
 \frac{d}{dt} (T_tf,\phi)=(LT_tf,\phi), \quad t\ge 0.
 \end{equation}
\end{theorem}

\Pf First let us calculate $\nabla_j \E g(X_x^{\tau}(k\tau))$ for an
arbitrary $k$ and $g\in C^1_{Lip}(\R^d)$. One has
\[
\E g(X_x^{\tau}(k\tau))=\int
g(x+z_1+...+z_k)P^{\tau}_x(dz_1)...P^{\tau}_{x+\sum_{m=1}^{k-1}z_m}(dz_k).
\]
As
\[
\nabla_j \E g(X_x^{\tau}(k\tau))= \lim_{h\to 0}
 \frac{1}{h}(\E g(X_{x+he_j}^{\tau}(k\tau))-\E g(X_x^{\tau}(k\tau)))
 \]
 does not depend on coupling, one can write
\[
\nabla_j \E g(X_x^{\tau}(k\tau))= \lim_{h\to 0}
 \frac{1}{h}\int (g(x+he_j+w_1+...+w_k)-g(x+v_1+...+v_k))
 \]
 \[
 P^{\tau}_{x+he_j,x}(dw_1dv_1)
...P^{\tau}_{x+he_j+\sum_{m=1}^{k-1}w_m,x+\sum_{m=1}^{k-1}v_m}(dw_kdv_k)
\]
\[
=\lim_{h\to 0}
 \int \bigl[(\nabla g(x+v_1+...+v_k),e_j+\frac{w_1-v_1}{h}...+\frac{w_k-v_k}{h})
 \]
 \[
 +O(1)\frac{1}{h} (he_j+w_1-v_1+...+w_k-v_k)^2\bigr]
  \]
 \[
 P^{\tau}_{x+he_j,x}(dw_1dv_1)
...P^{\tau}_{x+he_j+\sum_{m=1}^{k-1}w_m,x+\sum_{m=1}^{k-1}v_m}(dw_kdv_k).
\]
 The term with $O(1)$ vanishes as it can be rewritten by Proposition
\ref{proplevycoupling} as
\[
\lim_{h\to 0} O(1)\frac{1}{h}\int
(he_j+w_1-v_1+...+w_{k-1}-v_{k-1})^2(1+c\tau)
\]
\[
P^{\tau}_{x+he_j,x}(dw_1dv_1)
...P^{\tau}_{x+he_j+\sum_{m=1}^{k-2}w_m,x+\sum_{m=1}^{k-2}v_m}(dw_{k-1}dv_{k-1}),
\]
and consequently, iterating this procedure as
\[
\lim_{h\to 0} \frac{1}{h} O(1) h^2 (1+c\tau)^k=0.
\]
Hence
\[
\nabla_j \E g(X_x^{\tau}(k\tau)) =\lim_{h\to 0}
 \frac{1}{h}\int (\nabla g(x+v_1+...+v_k),he_j+w_1-v_1+...+w_k-v_k)
 \]
 \[
 P^{\tau}_{x+he_j,x}(dw_1dv_1)
...P^{\tau}_{x+he_j+\sum_{m=1}^{k-1}w_m,x+\sum_{m=1}^{k-1}v_m}(dw_kdv_k).
\]
Assume now first that $g$ is from the Schwartz space $S(\R^d)$ so
that Proposition \ref{proplevyderivcont1} applies and one can write
\[
\int (\nabla
g(x+v_1+...+v_k),he_j+w_1-v_1+...+w_k-v_k)P^{\tau}_{x+he_j+\sum_{m=1}^{k-1}w_m,x+\sum_{m=1}^{k-1}v_m}(dw_kdv_k)
\]
\[
=\int \sum_{j_k,j_{k-1}=1}^d \nabla_{j_k}
g(x+v_1+...+v_k)(he_j+w_1-v_1+...+w_{k-1}-v_{k-1})^{j_{k-1}}(\delta_{j_{k-1}}^{j_k}+z_k^{j_k})
\]
\[
Q^{\tau}_{D^{j_{k-1}}\nu(x+\sum_{m=1}^{k-1}v_m)}(dz_kdv_k)
+O(\tau)(w_{k-1}-v_{k-1})^2.
\]
Consequently, as the last term does not contribute to the limit
$h\to 0$, and iterating this procedure one obtains
\[
 \nabla_j \E
g(X_x^{\tau}(k\tau)) =
 \int \sum_{j_1,...,j_k=1}^d \nabla_{j_k} g(x+v_1+...+v_k)(\delta_{j_1}^j+z_1^{j_1})
  (\delta_{j_2}^{j_1}+z_2^{j_2})...(\delta_{j_{k-1}}^{j_k}+z_k^{j_k})
 \]
 \begin{equation}
\label{eqderivforfunctionsonapproximations}
 Q^{\tau}_{D^j\nu(x)}(dz_1dv_1)Q^{\tau}_{D^{j_1}\nu(x+v_1)}(dz_2dv_2)...
Q^{\tau}_{D^{j_{k-1}}\nu(x+\sum_{m=1}^{k-1}v_m)} (dz_k dv_k),
\end{equation}
which is the rigorous explicit form of the (a priori not clearly
defined but intuitively appealing) expression
\[
\E \nabla g(x+Y^0_{\tau}(x)+Y^1_{\tau}(X (\tau))+...+
Y^{k-1}_{\tau}(X ((k-1)\tau))
\]
\[
\left(1+\frac{\pa Y^{k-1}_{\tau}(X((k-1)\tau))}{\pa
X((k-1)\tau)}\right) ...\left(1+\frac{\pa Y^0_{\tau}(x)}{\pa
x}\right).
\]
 Approximating arbitrary $g$ by functions from the Schwartz space
one can conclude that \eqref{eqderivforfunctionsonapproximations}
holds for all $g\in C^1_{Lip}(\R^d)$.

 We want to show now that these derivatives are
Lipschitz continuous. To shorten the formulas let us do it for the
case of $d=1$ only. In this case
\[
 \nabla \E
g(X_x^{\tau}(k\tau)) =
 \int \nabla g(x+v_1+...+v_k)(1+z_1)...(1+z_k)
 \]
 \begin{equation}
\label{eqderivforfunctionsonapproximationsdim1}
 Q^{\tau}_{D\nu(x)}(dz_1dv_1)Q^{\tau}_{D\nu(x+v_1)}(dz_2dv_2)...
Q^{\tau}_{D\nu(x+\sum_{m=1}^{k-1}v_m)} (dz_k dv_k),
\end{equation}
and by Proposition \ref{proplevyderivcont} one can write
\[
\nabla \E g(X_{x_1}^{\tau}(k\tau))-\nabla \E
g(X_{x_2}^{\tau}(k\tau))
 \]
  \[
=\int \left[ \nabla g(x_1+v_1+...+v_k)(1+z_1)...(1+z_k)-\nabla
g(x_2+\tilde v_1+...+\tilde v_k)(1+\tilde z_1)...(1+\tilde
z_k)\right]
\]
\[
Q^{\tau}_{D\nu(x_1,x_2)}(dz_1d\tilde z_1dv_1d\tilde v_1)
 Q^{\tau}_{D\nu(x_1+v_1,x_2+\tilde v_1)}(dz_2d\tilde z_2dv_2d\tilde v_2)...
Q^{\tau}_{D\nu (x_1+\sum_{m=1}^{k-1}v_m, x_2+\sum_{m=1}^{k-1}\tilde
v_m)} (dz_kd\tilde z_k dv_kd \tilde v_k).
\]
Writing
\[
\nabla g(x_1+v_1+...+v_k)(1+z_1)...(1+z_k)-\nabla g(x_2+\tilde
v_1+...+\tilde v_k))(1+\tilde z_1)...(1+\tilde z_k)
\]
\[
=(\nabla g(x_1+v_1+...+v_k)-\nabla g(x_2+\tilde v_1+...+\tilde
v_k))(1+z_1)...(1+z_k)
 \]
 \[
 +\nabla g(x_2+\tilde v_1+...+\tilde v_k)[(1+z_1)...(1+z_k)-(1+\tilde
z_1)...(1+\tilde z_k)],
 \]
 and applying the H\"older inequality to estimate the integral over each of these two
 terms yields the estimate
\[
|\nabla \E g(X_{x_1}^{\tau}(k\tau))-\nabla \E
g(X_{x_2}^{\tau}(k\tau))| \le \kappa
[I_0^2(k,x_1,x_2)+I_1^2(k,x_1,x_2)]
 \]
 with $\kappa$ depending on the norm and the Lipschitz constant of
 $\nabla g$, where
\[
I_0^2(k,x_1,x_2)= \int [x_1-x_2+v_1-\tilde v_1+...+v_k-\tilde v_k]^2
|(1+z_1)^2...(1+z_k)^2
\]
\[
Q^{\tau}_{D\nu(x_1,x_2)}(dz_1d\tilde z_1dv_1d\tilde v_1)
 ...
Q^{\tau}_{D\nu (x_1+\sum_{m=1}^{k-1}v_m, x_2+\sum_{m=1}^{k-1}\tilde
v_m)} (dz_kd\tilde z_k dv_kd \tilde v_k).
\]
and
\[
I_1^2(k,x_1,x_2)= \int [(1+z_1)...(1+z_k)-(1+\tilde z_1)...(1+\tilde
z_k)]^2
\]
\[
Q^{\tau}_{D\nu(x_1,x_2)}(dz_1d\tilde z_1dv_1d\tilde v_1)
 ...
Q^{\tau}_{D\nu (x_1+\sum_{m=1}^{k-1}v_m, x_2+\sum_{m=1}^{k-1}\tilde
v_m)} (dz_kd\tilde z_k dv_kd \tilde v_k).
\]
By \eqref{eqlevyderivcont1}
\[
I_0^2(k,x_1,x_2) \le (1+c\tau)I_0^2(k-1,x_1,x_2) \le ... \le
(1+c\tau)^k(x_1-x_2)^2.
\]
It remains to estimate $I_1^2$. It would be convenient here to
introduce special notations for the products:
\[
Z_k=(1+z_1)...(1+z_k), \quad \tilde Z_k=(1+\tilde z_1)...(1+\tilde
z_k).
\]
Now one can write
\[
Z_k-\tilde Z_k=Z_{k-1}(1+z_k)-\tilde Z_{k-1}(1+\tilde z_k)
 =(1+z_k)(Z_{k-1}-\tilde Z_{k-1})+(z_k-\tilde z_k)\tilde Z_{k-1},
 \]
 so that
\[
[Z_k-\tilde Z_k]^2
 =(1+z_k)^2(Z_{k-1}-\tilde Z_{k-1})^2+(z_k-\tilde z_k)^2\tilde Z_{k-1}^2
 +2(1+z_k)(z_k-\tilde z_k)(Z_{k-1}-\tilde Z_{k-1})\tilde Z_{k-1}.
 \]
Plugging this into the expression for $I_1^2$ yields
\[
I_1^2(k,x_1,x_2) \le (1+c\tau)I_1^2(k-1,x_1,x_2)+c\tau
I_0^2(k-1,x_1,x_2)
 +c\tau
  \int \Omega (Z_{k-1}-\tilde Z_{k-1})\tilde Z_{k-1}
  \]
  \[
 Q^{\tau}_{D\nu(x_1,x_2)}(dz_1d\tilde z_1dv_1d\tilde v_1)
... Q^{\tau}_{D\nu (x_1+\sum_{m=1}^{k-2}v_m,
x_2+\sum_{m=1}^{k-2}\tilde v_m)} (dz_{k-1}d\tilde z_{k-1} dv_{k-1}d
\tilde v_{k-1}),
\]
where $\Omega$ in the last integral is a function of
$x_1,x_2,v_j,\tilde v_j$ such that
\[
|\Omega| \le c \|x_1-x_2+v_1-\tilde v_1+...+v_{k-1}-\tilde
v_{k-1}\|.
\]
Hence, applying to this last integral again the H\"older inequality
yields
\[
I_1^2(k,x_1,x_2) \le (1+c\tau)I_1^2(k-1,x_1,x_2)+c\tau
I_0^2(k-1,x_1,x_2),
 \]
 which taking into account the above bound for $I_0^2$ rewrites as
 \[
I_1^2(k,x_1,x_2) \le (1+c\tau)I_1^2(k-1,x_1,x_2)+c\tau (x_1-x_2)^2
 \]
 with yet another $c$ as long as $t=\tau k$ remains bounded.
 Using this formula recursively implies
 \[
I_1^2(k,x_1,x_2) \le c\tau (x_1-x_2)^2(1+(1+c\tau)+...+(1+c\tau)^k)
\le c(k\tau) (x_1-x_2)^2.
\]
Consequently one obtains the uniform estimate
\[
|\nabla \E g(X_{x_1}^{\tau}(k\tau))-\nabla \E
g(X_{x_2}^{\tau}(k\tau))| \le \kappa c(k\tau)\|x_1-x_2\|.
 \]
 Hence from the sequence of the uniformly Lipschitz continuous functions
$\nabla \E f(X^{\tau_k}_x(s))$, $k=1,2,...$, one can choose a
convergent subsequence the limit being clearly $\nabla \E
f(X_x(t))$, showing that $\E f(X_x(t))\in C^1_{Lip}$. The uniform
convergence implies $\E f(X_x(t))\in C^1_{Lip}\cap C^1_{\infty}$
whenever the same holds for $f$.

To complete the proof of the theorem it remains equation
\eqref{eqevolutioneqforsemigroupinLlip}. But this is easy: for $t=0$
it follows by  approximating $f$ with $f_n\in C^2(\R^d)$ and then
for arbitrary $t$ it follows by the invariance of the class
$C^1_{Lip}$ under $T_t$.

 Second derivative can be analyzed similarly, but the assumptions and calculations become essentially longer.

 \appendix

 \section{Coupling of L\'evy processes}
 \label{appendix2}

 We describe here the natural coupling of L\'evy processes leading in particular to the
 analysis
of their weak derivatives with respect to a parameter. Recall that
by $C^k_{Lip}$ we denote the subspace of functions from $C^k(\R^d)$
with a Lipschitz continuous derivative of order $k$.

\begin{prop}
\label{proplevycoupling}

Let $Y_s^i$, $i=1,2$, be two L\'evy processes in $\R^d$ specified by
their generators
\begin{equation}
\label{levygenerator}
 L_if(x)=\frac{1}{2}(G_i\nabla,\nabla)f(x)+
(b_i,\nabla f(x))
 +\int (f(x+y)-f(x)-(\nabla f (x), y))\nu_i (dy)
 \end{equation}
 with $\nu_i \in \MC_2(\R^d)$. Let $\nu \in \MC_2(\R^{2d})$ be a coupling of $\nu_1,\nu_2$, i.e.
\eqref{eqcouplingcond} holds for all $\phi_1,\phi_2$ satisfying
$\phi_i(.)/|.|^2 \in C(\R^d)$. Then the operator
\[
 Lf(x_1,x_2)=\left[\frac{1}{2}(G_1\nabla_1,\nabla_1)
 + \frac{1}{2}(G_2\nabla_2,\nabla_2)
  +(\sqrt {G_2} \sqrt {G_1} \nabla_1,\nabla_2)
  \right]f(x_1,x_2)
\]
\[
 + (b_1,\nabla_1 f(x_1,x_2)) + (b_2,\nabla_2 f(x_1,x_2))
\]
\begin{equation}
\label{levycouplinggenerator2}
 +\int [f(x_1+y_1,x_2+y_2)-f(x_1,x_2)-((y_1,\nabla_1)+(y_2,\nabla_2))f(x_1,x_2)]
 \nu(dy_1dy_2)
 \end{equation}
 (where $\nabla_i$ means the gradient with respect to $x_i$)
specifies a L\'evy process $Y_s$ in $\R^{2d}$ with the
characteristic exponent
\[
 \eta_{x_1,x_2}(p_1,p_2)
 =-\frac{1}{2} \left[
 \sqrt {G(x_1)}p_1+\sqrt {G(x_2)}p_2\right]^2
+ib(x_1)p_1+ib(x_2)p_2
\]
\[
 +\int (e^{iy_1p_1+iy_2p_2}-1-i(y_1p_1+y_2p_2))
 \nu(dy_1dy_2),
\]
that is a coupling of $Y^1_s,Y^2_s$ in the sense that the components
of $Y_s$ have the distribution of $Y^1_s$ and $Y^2_s$ respectively.
 Moreover, if $f(x_1,x_2)=h(x_1-x_2)$ with a function $h\in C^2(\R^d)$,
 then
\[
 Lf(x_1,x_2)=\frac{1}{2}((\sqrt {G_1}-\sqrt
 {G_2})^2\nabla,\nabla)h(x_1-x_2)+(b_1-b_2,\nabla h)(x_1-x_2)
\]
\begin{equation}
\label{levycouplinggenerator2a}
 +\int [h(x_1-x_2+y_1-y_2)-h(x_1-x_2)-(y_1-y_2,\nabla h(x_1-x_2)]
 \nu(dy_1dy_2).
 \end{equation}
 In particular, if $f(x_1,x_2)=(x_1-x_2)^2$, then
 \begin{equation}
\label{levycouplinggenerator2b}
 Lf(x_1,x_2)=\tr (\sqrt {G_1}-\sqrt
 {G_2})^2+2(b_1-b_2,x_1-x_2)
 +\int (y_1-y_2)^2\nu(dy_1dy_2).
 \end{equation}
  Finally
\begin{equation}
\label{levycouplingestimate2}
 \E (\xi+Y_t^1-Y_t^2)^2
=(\xi+t(b_1-b_2))^2
 +t\left(Tr (\sqrt{G_1}-\sqrt {G_2})^2
 +\int \int (y_1-y_2)^2\nu(dy_1dy_2)\right).
 \end{equation}
 \end{prop}

 \Pf Straightforward. In fact, clearly $Y_s$ couples
 $Y_s^1,Y_s^2$, because say $\eta_{x_1,x_2}(p_1,0)$ is the characteristic exponent
of $Y_s^1$. Equation \eqref{levycouplinggenerator2a} follows from
\eqref{levycouplinggenerator2}. The second moment
\eqref{levycouplingestimate2} is found either by twice
differentiating the characteristic function, or by the Dynkin
formula in conjunction with \eqref{levycouplinggenerator2a}.

Similarly one obtains

\begin{prop}
\label{proplevycoupling1}

Let $Y_s^i$, $i=1,2$, be two L\'evy processes in $\R^d$ specified by
their generators
\begin{equation}
\label{levygenerator1}
 L_if(x)= (b_i,\nabla f(x)) +\int
(f(x+y)-f(x))\nu_i (dy)
 \end{equation}
 with $\nu_i \in \MC_1(\R^d)$. Let $\nu \in \MC_1(\R^{2d})$ be a coupling of $\nu_1,\nu_2$, i.e.
\eqref{eqcouplingcond} holds for all $\phi_1,\phi_2$ satisfying
$\phi_i(.)/|.| \in C(\R^d)$. Then the operator
\begin{equation}
\label{levycouplinggenerator1}
 Lf(x_1,x_2)=(b_1,\nabla_1 f(x_1,x_2))
+ (b_2,\nabla_2 f(x_1,x_2))
 +\int [f(x_1+y_1,x_2+y_2)-f(x_1,x_2)] \nu(dy_1dy_2)
 \end{equation}
specifies a L\'evy process $Y_s$ in $\R^{2d}$ that is a coupling of
$Y^1_s,Y^2_s$ such that for all $t$
\begin{equation}
\label{levycouplingestimate1}
 {\bf E}\|\xi+Y_t^1-Y_t^2\| \le \|\xi \|
 +t\left(\|b_1-b_2\|
 +\int \int \|y_1-y_2\|\nu(dy_1dy_2)\right).
 \end{equation}
 \end{prop}

\Pf One approximates $|y|$ by a smooth function, applies Dynkin's
formula and then passes to the limit.

Next, let $Y_t(z)$ be a family of L\'evy processes in $\R^d$
parametrized by points $z\in \R^d$ and specified by their generators
\begin{equation}
\label{eqfamilyoflevygenerators}
L[z]f(x)=\frac{1}{2}(G(z)\nabla,\nabla)f(x)+ (b(z),\nabla f(x))
 +\int (f(x+y)-f(x)-(\nabla f (x), y))\nu (z;dy)
 \end{equation}
 where $\nu (z;.) \in \MC_2(\R^d)$.
We are interested in defining the process $\frac{\pa}{\pa z}
Y_t(z)$.

We shall describe this process via a certain derivative type
operator on L\'evy measures connected with a coupling. Namely, let
$\nu_{x_1,x_2}(dy_1dy_2)$ be a family of $\MC_2$-couplings of
$\nu(x_1;.),\nu(x_2;,.)$ (in the sense that $\nu_{x_1,x_2} \in
\MC_2(\R^{2d})$ and \eqref{eqcouplingcond} holds for all
$\phi_1,\phi_2$ such that $\phi_i(.)/|.|^2 \in C(\R^d)$). For
instance, these could be optimal couplings with respect to the cost
function $(y_1-y_2)^2$, i.e. those couplings, where the infinum in
the definition of $W_2(\nu(x_1,.),\nu(x_2,.))$ is attained.

Let $T_h(y_1,y_2)=((y_1-y_2)/h,y_2)$ and the measure
$\nu^{T_h}_{T_h^{-1}(\xi,x)}$ on $\R^{2d}$ be defined as the push
forward of $\nu_{x+h\xi,x}=\nu_{T_h^{-1}(\xi,x)}$ by $T_h$, i.e.
\[
\int \int f(z,y)\nu^{T_h}_{T_h^{-1}(\xi,x)} (dz dy)= \int \int
f(\frac{y_1-y_2}{h},y_2)\nu_{x+h\xi,x}(dy_1dy_2).
\]
Clearly $\nu^{T_h}_{T_h^{-1}(\xi,x)}$ is a L\'evy measure with a
finite second moment whenever this is the case for $\nu_{x+h\xi,x}$.
The relevant smoothness of $\nu$ will be defined now as the
existence of the weak limit
\[
D_{\xi}\nu_x =\lim_{h\to 0}\nu^{T_h}_{T_h^{-1}(\xi,x)},
\]
i.e.
\[
\lim_{h\to 0} \int \int g(y_1,y_2)
\nu^{T_h}_{T_h^{-1}(\xi,x)}(dy_1dy_2)=\int \int
g(y_1,y_2)D_{\xi}\nu_x (dy_1dy_2), \quad \frac
{g(y_1,y_1)}{y_1^2+y_2^2} \in C(\R^{2d}).
\]
To see the rational behind this definition observe that if
$\nu(x,.)=\nu^{F_x}$ with a given $\nu$ and a family of
transformations $F_x(.)$, then
 \[
D_{\xi}\nu_x=\nu^{(\xi,\nabla F_x(.)),F_x(.)}
\]
is the push forward of $\nu$ with respect to $y\mapsto
 (\xi,\nabla F_x(y)),F_x(y)$ ($\nabla$ is the derivative with respect to $x$).
  On the other hand, if $\nu_{z,x}$ has a density,
 i.e. $\nu_{z,x}(dy_1dy_2)=\nu_{z,x}(y_1,y_2)dy_1dy_2$. then
$D_{\xi}\nu_x$ has the density
\[
\lim_{h\to 0} h^d \nu_{x+h\xi,x}(y+hz,y).
\]
If the coupling is optimal (is given by minimizers in the definition
of the $W_2$- distance) this derivative is connected with the
derivative of $W_2$
 via the formula
 \[
 \int |z|^2
D_{\xi}\nu_x (dzdy)=\left( \frac{d}{dh}\mid_{h=0} W_2(\nu(x+h\xi;.),
\nu(x;.))\right)^2.
\]

We shall need further only the partial derivatives
$D^i\nu_x=D_{e_i}\nu_x$ in the directions of the co-ordinate vectors
$e_i$. The reason for introducing these derivatives lies in the
observation that its action on L\'evy measures corresponds to the
derivation of L\'evy processes. More precisely, the following holds.

\begin{prop}
\label{proplevyderiv} Let $Y_t(z)$ be the family of the L\'evy
processes in $\R^d$, $z\in \R^d$, specified by their generators
\eqref{eqfamilyoflevygenerators}. Suppose $G(x),b(x) \in C^1(\R^d)$
and $\nu(x,.)$ is smooth in the above sense (i.e. $ D^j\nu$ are well
defined with respect to a certain coupling). (i) Then the coupled
random variables in $\R^{2d}$
\[
(h^{-1}(Y_t(x+he_j)-Y_t(x)), Y_t(x))
\]
in $\R^{2d}$ has a weak limit that we denote
$(\nabla_jY_t(x),Y_t(x))$ and that has the distribution
$Q^t_{D^j\nu(x)} $ of the L\'evy process at time $t$ with the
characteristic exponent
\[
 \eta_x^j(q,p)=-\frac{1}{2} \left[\nabla_j\sqrt {G(x)}q+\sqrt {G(x)}p\right]^2
+i(\nabla_jb(x),q)+i(b(x),p)
\]
\begin{equation}
\label{eqcharexplevyderivative}
 +\int (e^{iqz+ipy}-1-ipy-iqz)) D^j\nu_x (dzdy).
\end{equation}
(ii) Moreover, if $g\in C^1_{Lip} (\R^{2d})$, then the partial
derivatives $\nabla_j \E g(x,Y_t(x))$ exist and
\begin{equation}
\label{eqlevyderivativeparameter}
 \nabla_j \E g(x,Y_t(x)) =\int
\left(\nabla_j g (x,y) +(\frac{\pa g}{\pa y}
(x,y),z)\right)Q^t_{D^j\nu(x)} (dz dy)
\end{equation}
($\nabla_j$ means the derivative with respect to the variable $x$).
\end{prop}

\Pf (i) The characteristic exponent of the L\'evy process
$T_h(Y_t(T_h^{-1}(e_j,x))$ is
\[
 \eta_x^{j,h}(q,p)=-\frac{1}{2} \left[\sqrt {G(x+he_j)}\frac{q}{h}+\sqrt {G(x)}(p-\frac{q}{h})\right]^2
+i(b(x+he_j),\frac{q}{h})+i(b(x),(p-\frac{q}{h})
\]
\[
 +\int (e^{iy_1q/h+iy_2(p-q/h)}-1-i(y_1-y_2)\frac{q}{h}-ipy_2))
 \nu_{T_h^{-1}(e_j,x)}(dy_1dy_2),
\]
which clearly converges to \eqref{eqcharexplevyderivative}.

 (ii) One has
 \[
 \frac{1}{h} \left[\E g(x+he_j,Y_t(x+he_j))-\E g(x,Y_t(x))\right]
 =\frac{1}{h} \E \left[g(x+he_j,Y_t(x+he_j))-g(x,Y_t(x))\right],
 \]
 where the last expectation can be taken with respect to any coupling
 of $Y_t(x+he_j)$ and $Y_t(x)$. Hence it can be written as
 \[
 \int\int \left(\nabla_j g (x,y_2)+(\frac{\pa g}{\pa y} (x,y_2),\frac{y_1-y_2}{h})
 +O(1) \frac{1}{h} (h^2\xi^2 +(y_1-y_2)^2\right)P^t_{x+he_j,x}(dy_1dy_2).
 \]
 By the property of the coupling (Proposition \ref{proplevycoupling}) the term with $O(1)$ tends to zero
 as $h\to 0$. Consequently
 \[
\frac{d}{dh}\mid_{h=0}\E g(x+he_j,Y_t(x+he_j))
 =\int \nabla_j g(x,y)P^t_x(dy)
 \]
 \[
 +\lim_{h\to 0}
\int \int (\frac{\pa g}{\pa y} (x,y_2),\frac{y_1-y_2}{h})
P^t_{x+he_j,x}(dy_1dy_2),
\]
implying \eqref{eqlevyderivativeparameter} due to statement (i).

It is worth noting that statement (ii) implies that the
distributions of the derivatives actually do not depend on coupling.

So far we have got only partial derivatives. We are now interested
in their continuity which clearly is linked to the continuity of the
measures $D^i\nu_x$. It turns out that the relevant notion of
continuity is a bit finer than the $W_2$-continuity used above. Next
two statements reveal two 'crucial bits' of this continuity.

\begin{prop}
\label{proplevyderivcont} Under the assumptions of Proposition
\ref{proplevyderiv} assume additionally that $G(x),b(x) \in
C^1_{Lip}(\R^d)$ and that the L\'evy measures $D^j\nu_x$ are
Lipschitz continuous in the following sense: for any $x_1,x_2\in
\R^d$ and $j=1,...,d$ there exists  a L\'evy coupling $D^j(x_1,x_2)$
of the L\'evy measures  $D^j\nu_{x_1}$, $D^j\nu_{x_2}$ such that
\begin{equation}
\label{eqlevyderivcont}
 \int_{\R^{4d}} \left[
(y_1-y_2)^2(1+z_1^2+z_2^2)+(z_1-z_2)^2\right]
D^j(x_1,x_2)(dz_1dz_2dy_1dy_2) \le \kappa (x_1-x_2)^2
\end{equation}
with a constant $\kappa$. Let $Q^t_{D^j(x_1,x_2)}$ denote the
distribution at time $t$ of the L\'evy process that couples
$(\nabla_jY_t(x_1),Y_t(x_1))$ and $(\nabla_jY_t(x_2),Y_t(x_2))$
according to Proposition \ref{proplevycoupling}, i.e. the L\'evy
process in $\R^{4d}$ specified by the characteristic exponent
\[
 \eta^j_{x_1,x_2}(q_1,q_2,p_1,p_2)
 =-\frac{1}{2} \left[\nabla_j\sqrt {G(x_1)}q_1+\nabla_j\sqrt
 {G(x_2)}q_2
 +\sqrt {G(x_1)}p_1+\sqrt {G(x_2)}p_2\right]^2
\]
\[
+i(\nabla_jb(x_1)q_1+\nabla_jb(x_2)q_2+b(x_1)p_1+b(x_2)p_2)
\]
\[
 +\int (e^{iy_1p_1+iy_2p_2+iz_1q_1+iz_2q_2}-1-i(y_1p_1+y_2p_2+z_1q_1+z_2q_2)
 D^j_{x_1,x_2}(dz_1dz_2dy_1dy_2).
\]
Then for any $\xi \in \R^d$
\begin{equation}
\label{eqlevyderivcont1}
 \int_{\R^{4d}} \left[
(\xi +y_1-y_2)^2(1+z_1)^2+(z_1-z_2)^2\right]
Q^t_{D^j(x_1,x_2)}(dz_1dz_2dy_1dy_2) \le \xi^2+ct(\xi^2+(x_1-x_2)^2)
\end{equation}
with a constant $c$ uniformly for finite times, and for any $g\in
C^1_{Lip}(\R^{2d})$ the function $\E g(x,Y_t(x))$ belongs to
$C^1_{Lip}(\R^d)$ (also uniformly for finite times).
\end{prop}

\Pf The moment estimates \eqref{eqlevyderivcont1} are obtained
directly from the derivatives of the characteristic function as in
Proposition \ref{proplevycoupling}. For the second statement we
write
\[
|\nabla_j\E g(x_1,Y_t(x_1))-\nabla_j\E g(x_1,Y_t(x_1))|
 \le \int_{\R^{4d}} Q^t_{D^j(x_1,x_2)}(dz_1dz_2dy_1dy_2)
 \]
 \[
\times |\nabla_j g(x_1,y_1)-\nabla_j g(x_2,y_2)+(\frac{\pa g}{\pa
y}(x_1,y_1),z_1)+(\frac{\pa g}{\pa
 y}(x_2,y_2),z_2)|
 \]
 (the derivative $\nabla_j$ with respect to $x$),
 which does not exceed
\[
\int_{\R^{4d}} ((|x_1-x_2|+|y_1-y_1|(1+|z_1|+|z_2|)+|z_1-z_2|)
Q^t_{D^j(x_1,x_2)}(dz_1dz_2dy_1dy_2),
 \]
 and which in turn does not exceed $\sqrt t |x_1-x_2|$ due to \eqref{eqlevyderivcont1}
 and the H\"older inequality.

In case $\nu(x,.)=\nu^{F_x}$ for a family of transformations
$F_x(.)$ the coupling $D^j(x_1,x_2)$ can be obtained as
\[
\int f(z_1,z_2,y_1,y_2) D^j(x_1,x_2)(dz_1dz_2dy_1dy_2)
 \]
  \[
 =\int f
(\nabla_jF(x_1,y), \nabla_j F(x_2,y), F(x_1,y), F(x_2,y)) \nu (dy),
\]
and the condition \eqref{eqlevyderivcont} is fulfilled whenever the
derivatives $\frac{\pa} {\pa x}F(x,y)$ are bounded and Lipschitz
continuous.

By $D\nu_x$ we shall denote the vector $\{D^j\nu_x\}$ and by
$Q^t_{D\nu(x)}$ the vector $\{Q^t_{D^j\nu (x)}\}$, $j=1,...,d$.

\begin{prop}
\label{proplevyderivcont1}
 Under the assumptions of Proposition
\ref{proplevyderiv} assume additionally that $G(x),b(x) \in
C^1_{Lip}(\R^d)$ and that the function
\[
 \int\int (y_1-y_2,e_j)(e^{iy_2p}-1)
\nu_{x,z}(dy_1dy_2)
\]
 is differentiable in $x$ around $x=z$ with
uniform estimates, more precisely that
\[
 \int\int (y_1-y_2,e_j) (e^{iy_2p}-1)\nu_{x,z}(dy_1dy_2)
\]
\begin{equation}
\label{eqlevyderivcont2}
 =\left(\frac{\pa}{\pa x}\mid_{x=z} \int\int (y_1-y_2,e_j) (e^{iy_2p}-1)
\nu_{x,z}(dy_1dy_2),x-z\right)+O(1+|p|)(x-z)^2.
\end{equation}
 Then for a continuous function $g$ represented via the inverse
Fourier transform as
\[
g(y)=\int e^{iyp} \hat g(p) \, dp, \quad (1+|p|)\hat g (p) \in
(L^1(\R^d))^d,
\]
one has the estimate
\[
 \E (Y_t(x)-Y_t(z), g(Y_t(z)))=\int
(y_1-y_2, g(y_2)) P^{\tau}_{x,z}(dy_1dy_2)
\]
\begin{equation}
\label{eqlevyderivcont3}
 =\int\int (w,g(y)) (Q^t_{D\nu (z)} (dw dy), x-z) +O(t)(x-z)^2\int (1+|p|)|\hat g (p)| \, dp.
\end{equation}
\end{prop}

\Pf Comparing the r.h.s of \eqref{eqlevyderivcont2} with the
definition of $D\nu_x$ yields
\[
 \int\int (y_1-y_2,e_j) (e^{iy_2p}-1)\nu_{x,z}(dy_1dy_2)
 \]
\begin{equation}
\label{eqlevyderivcont5}
 = \left(\int\int (w,e_j) (e^{iyp}-1)
D\nu_x(dw dy),x-z\right) +O(1+|p|)(x-z)^2.
\end{equation}
Now one has
\[
\int (y_1-y_2,e_j)e^{iy_2p} P^t_{x,z}(dy_1dy_2)
\]
\[
 =-i \frac{\pa}{\pa q^j}\mid_{q=0}\E \exp
 \{i(Y_t(x)-Y_t(z))q+iY_t(z)p\}=
-i \frac{\pa}{\pa q^j}\mid_{q=0}\exp \{t\eta_{x,z}(q,p-q)\}
 \]
 \[
  =t\left[i (\sqrt {G(z)} (\sqrt {G(x)}-\sqrt {G(z)})p)^j+ (b(x)-b(z))^j +\int
(y_1-y_2)^j(e^{ipy_2}-1)\nu_{x,z}(dy_1dy_2)\right]
 \E e^{iY_t(z)p}
\]
\[
=t\left(\frac{i}{2}(\nabla (G(z)p)^j+ \nabla b^j(z) +\int\int w^j
(e^{ipy}-1) (D\nu_x(dw dy),x-z\right)
 \E e^{iY_t(z)p}+O(t)(1+|p|)(x-z)^2.
 \]
 Consequently,
\[
\int (y_1-y_2,g(y_2)) P^t_{x,z}(dy_1dy_2)
\]
\[
=t\int \left(\frac{i}{2}\nabla (G(z)p, \hat g(p))+ \nabla (b(z),
\hat g(p)) +\int\int (w, \hat g(p)) (e^{ipy}-1) D_{\nu_x}(dw
dy),x-z\right)
 \E e^{iY_t(z)p} \, dp
 \]
 \[
  +O(t)\int (1+|p|) \hat g(p) \, dp (x-z)^2.
 \]
 Similarly
 \[
 \int \int w^j e^{ipy} Q^t_{D\nu (z)} (dw dy)=
 -i\frac{\pa}{\pa q^j}\mid_{q=0}
 \]
 \[
 \exp \{t [-\frac{1}{2} (\nabla \sqrt {G(z)}q+\sqrt {G(z)} p]^2
 +i(\nabla b(z)q+b(z)+\int \int (e^{iqw+ipy}-1-ipy-iqw)D\nu_x(dw dy)\}
 \]
 \[
 = t \left(\frac{i}{2}\nabla (G(z)p)^j+ \nabla b^j(z) +\int\int w^j
(e^{ipy}-1) D\nu_x(dw dy)\right) \E e^{iY_t(z)p},
\]
 implying \eqref{eqlevyderivcont3}.

 To differentiate the L\'evy process for the second time, one needs
 of course the 'second derivative' of the L\'evy measure defined similarly to the first one. Namely,
one needs the existence of the limit
\begin{equation}
\label{eqsecondderLevymes}
 \lim_{h\to 0} \int f(\frac{z_1-z_2}{h},
\frac{y_1-y_2}{h},z_2,y_2)D^j(x+he_k,x)(dz_2dz_1dy_2dy_1) =\int
f(w,z_k,z_j,y)D^{kj}_x(dwdz_jdz_kdy)
\end{equation}
whenever $f(w,z_k,z_j,y)/(w^2+z_j^2+z_k^2+y^2) \in C(\R^{4d})$ with
$D^{kj}_x(dwdz_jdz_kdy)$ belonging to $\MC_2(\R^{4d})$. The
following is a straightforward analog of Proposition
\ref{proplevyderiv}.

\begin{prop}
\label{proplevysecondderiv} Under the assumptions of Proposition
\ref{proplevyderiv} assume that $G(x),b(x) \in C^2(\R^d)$ and the
measures $D^{kj}_x\in \MC_2(\R^{4d})$ are well defined by
\eqref{eqsecondderLevymes}. (i) Then for any $j,k$ the process
\[
(\nabla_k\nabla_j Y_t(x), \nabla_k Y_t(x), \nabla_j Y_t(x), Y_t(x))
\]
is defined weakly in $\R^{4d}$ and has the distribution
$Q^t_{D^{jk}\nu(x)}$ of the L\'evy process at time $t$ with the
characteristic exponent
\[
 \eta_x^{jk}(r,q_k,q_j,p)=-\frac{1}{2} \left[\nabla_k\nabla_j\sqrt {G(x)}q
+\nabla_k \sqrt {G(x)}q_k +\nabla_j \sqrt {G(x)}q_j
 +\sqrt {G(x)}p\right]^2
 \]
 \[
+i[(\nabla_k \nabla_jb(x),q)+(\nabla_kb(x),q_k)+(\nabla_j
b(x),q_j)+(b(x)p)]
\]
\begin{equation}
\label{eqcharexplevysecondderivative}
 +\int [e^{irw+iq_kz_k+iq_jz_j+ipy}-1-i(rw+q_kz_k+q_jz_j+py)] D^{kj}_x (dwdz_kdz_jdy).
\end{equation}
(ii) Moreover, if $g\in C^2_{Lip} (\R^d)$, the partial derivatives
$\nabla_k\nabla_j \E g(x+Y_t(x))$ exist and
\[
 \nabla_k\nabla_j \E g(x+Y_t(x)) =\int Q^t_{D^{kj}\nu(x)} (dwdz_kdz_jdy)
\]
\begin{equation}
\label{eqlevysecondderivativeparameter}
 \left[\sum_{l,m=1}^d \nabla_m\nabla_l g(x+y)(\delta_k^m+z_k^m)(\delta_j^l+z_j^l)
 +\sum_{l=1}^d \nabla g(x+y)w^l\right]
\end{equation}
($\nabla$ means the derivative with respect to the variable $x$).
\end{prop}

\paragraph{Acknowledgments.} The author is grateful to the referee
whose constructive critique of the first draft helped to improve
essentially the quality of the exposition.

\end{document}